\documentclass[a4paper,12pt]{article}
\usepackage{dsfont}
\usepackage{amsmath}

\usepackage{amssymb,amsthm}
\usepackage{mathrsfs,amsfonts}
\usepackage{color}
\usepackage{graphicx}
\usepackage{bbm}
\usepackage{ulem}

\allowdisplaybreaks

\allowdisplaybreaks

\newtheorem{con0}{Theorem}[section]
\newtheorem{thm0}{Theorem}[section]
\newtheorem{exa0}{Theorem}[section]

\newtheorem{con1}[con0]{Condition}
\newtheorem{def1}[thm0]{Definition}
\newtheorem{lem1}[thm0]{Lemma}
\newtheorem{thm1}[thm0]{Theorem}
\newtheorem{cor1}[thm0]{Corollary}
\newtheorem{pro1}[thm0]{Proposition}
\newtheorem{rem1}[thm0]{Remark}
\newtheorem{ass1}[thm0]{Assumption}
\newtheorem{exa1}[exa0]{\it{Example}}

\def\bglemma{\begin{lem1}}\def\edlemma{\end{lem1}}
\def\bgtheorem{\begin{thm1}}\def\edtheorem{\end{thm1}}

\def\bgproposition{\begin{pro1}}\def\edproposition{\end{pro1}}

\def\bgcondition{\begin{con1}}\def\edcondition{\end{con1}}

\def\benumerate{\begin{enumerate}}\def\eenumerate{\end{enumerate}}
\def\bitemize{\begin{itemize}}\def\eitemize{\end{itemize}}

\def\beqlb{\begin{eqnarray}}\def\eeqlb{\end{eqnarray}}
\def\beqnn{\begin{eqnarray*}}\def\eeqnn{\end{eqnarray*}}

\def\eqref#1{{\rm(\ref{#1})}}

\def\ar{\!\!\!&}\def\nnm{\nonumber}\def\ccr{\nnm\\}

\def\proof{\noindent{\it Proof.~}}
\def\qed{\hfill$\square$\smallskip}

\def\mrm{\mathrm}\def\mbb{\mathbf}\def\mcr{\mathscr}
\def\mbb{\mathbb}\def\mds{\mathds}

\def\d{\mrm{d}}\def\e{\mrm{e}}

\def\I{\mds{1}}

\def\itDelta{{\it\Delta}}\def\itPhi{{\it\Phi}}
\def\itPsi{{\it\Psi}}

\newcommand{\R}{\mathbb R}
\newcommand{\Pp}{\mathbb P}

\def\red{\color{red}}

\begin{document}

	\centerline{\Large\bf Quasi-stationary distribution for continuous-state}
	
	\smallskip
	
	\centerline{\Large\bf  branching processes with competition}
	
	\bigskip
	
	\centerline{Pei-Sen Li$^1$,  Jian Wang$^2$ and Xiaowen Zhou$^3$}
	
	\medskip
	
	{\small\it
		
		\centerline{$^1$School of Mathematics and Statistics, Beijing Institute of Technology, Beijing}
		
		\centerline{100872, China, {\tt peisenli@bit.edu.cn}}
		
		\smallskip

		\centerline{$^2$School of Mathematics and Statistics, Fuzhou 350007, China, {\tt jianwang@fjnu.edu.cn}}
		
		\smallskip
		
		\centerline{$^3$Department of Mathematics and Statistics, Concordia University, 1455 De}
		
		\centerline{Maisonneuve Blvd. W., Montreal, Canada, {\tt xiaowen.zhou@concordia.ca}}
		
	}
	
	\bigskip
	
	{\narrower{\narrower
			
			\noindent{\textit{Abstract:}} We study quasi-stationary distribution of the continuous-state branching process with competition introduced in Berestycki,  Fittipaldi and Fontbona\ (Probab. Theory Relat. Fields, 2018).  This process is constructed as the unique strong solution to a stochastic integral equation with jumps. An important example is the logistic branching process constructed in Lambert (Ann. Appl. Probab., 2005). We establish  the strong Feller property,
			trajectory Feller property, Lyapunov condition, weak Feller property and irreducibility, respectively. These properties together allow us to prove that if the competition  is strong enough near $+\infty$, then there is a unique quasi-stationary distribution, which attracts all initial  distributions  with exponential rates.
			
			\bigskip
			
			\noindent{\textit{Key words:}} continuous-state branching process;  competition;
			strong Feller property;
			irreducibility; quasi-stationary distribution.
			
			\smallskip
			
			\noindent{\textit{MSC {\rm(2020)} Subject Classification:}} 60J80, 60J25, 60G51, 60G52
			
			\par}\par}

	\bigskip

	\section{Introduction}
	\subsection{Background and main result}
	\setcounter{equation}{0}
	
	The main objective of this work is to establish the existence, uniqueness and exponential convergence of quasi-stationary distribution for a population model,  the continuous-state branching process with competition,  introduced by Berestycki et al.\ \cite{BFF18}. This process is constructed as the unique strong solution to a stochastic integral equation with jumps, which will be presented later. The quasi-stationary distribution is a crucial concept in population models.
It remains unchanged over time when the process is conditioned to survive and describes the asymptotic behavior of the process near extinction. For general introductions on this topic, we refer to the survey by M\'el\'eard and Villemonais \cite{MV12} and the  monograph by Collet et al. \cite{CMS}.
	
	Consider a non-negative Markov process $(Y_t)_{t\ge 0}$ on $E$ such that either $E=\mbb{N}=\{0, 1, 2,\dots\}$ or $E=\mbb{R}_+=[0, \infty)$.
		Assume that the state $0$ is an absorbing state; that is, $Y_t=0$ for $t\ge \tau^-_0$, where $\tau^-_0:=\inf\{t\ge0: Y_t=0\}$ is the extinction time with the convention $\inf\emptyset=\infty$.
	We say that a probability measure {$\pi$} on $E\backslash \{0\}$ is a quasi-stationary distribution (QSD for short) of $(Y_t)_{t\ge0}$ if, for all $t\ge0$ and any measurable set $A\subset E\backslash \{0\},$
	\[
	{\pi}(A)=\mbb{P}_\pi(Y_t\in A|\tau^-_0>t),
	\]
	where $\mbb{P}_\pi$ denotes the distribution of the process $Y$ with initial distribution $\pi.$
	
	From the definition above and the Markov property, we see that
	$\pi$ is a QSD if and only if it is a common, positive left-eigenvector of the killed transition semigroup; see \cite[Remark 2.3]{GNW}. So, in the finite state space, thanks to the Perrron-Frobenius theorem, we know there is a unique QSD; see \cite[Section 1.2]{CMMS}. However, for Markov process taking values in an infinite state space, even simple models such as branching processes may have infinitely many QSDs. Indeed, the study of QSD began with Yaglom's work \cite{Yaglom} on Bienaym\'e--Galton--Watson branching processes, which can be represented as
	\beqnn
	X_{n} = \sum_{i=1}^{X_{n-1} }\xi_{n,i}, \quad n\ge 1
	\eeqnn
 with $(\xi_{n,i})$
being i.i.d. random variables.
	Let $\rho=\mbb{E}(\xi_{n, i})$ be the mean of the offspring distribution. It was proved in Seneta and Vere-Jones \cite{SeV66} that  in the sub-critical situation, i.e., $\rho<1$,
	a Bienaym\'e--Galton--Watson branching process admits infinitely many QSDs. Later, Hope \cite{Hoppe} solved the problem of characterizing all QSDs; see also Maillard \cite{Mail18} for an analytical  approach to this problem.

	Similar to the relation between random walk and Brownian motion, it is natural
	to consider branching models with continuous-state space $\mbb{R}_+$. The study of \textit{continuous-state branching processes} (CB-processes for short) was initiated by Feller \cite{Fel51}, who noticed that a one-dimensional diffusion process may arise as the limit of a sequence of rescaled Bienaym\'e--Galton--Watson branching processes. The result was extended by Lamperti \cite{Lam67a} to the situation where the limiting process may have discontinuous sample paths. Let $\itPsi$ be a function on $[0,\infty)$ with the L\'evy--Khintchine representation:
	\beqlb\label{bran-mech}
	\itPsi(\lambda)= b\lambda + c\lambda^2 + \int_0^\infty \big(\e^{-\lambda z}-1+\lambda z\I_{\{z\le 1\}}\big)\mu(\d z), \quad \lambda\ge 0,
	\eeqlb
	where $b\in \mbb{R}$ and $c\ge 0$ are constants, and $(1\land z^2)\mu(\d z)$ is a finite measure on $(0,\infty)$. The transition kernel $Q_t(x,\d y)$ of a CB-process with \textit{branching mechanism} $\itPsi$ is defined by
	\beqlb\label{CB-semigr}
	\int_{\mbb{R}_+} \e^{-\lambda y} Q_t(x,\d y)= \e^{-xv_t(\lambda)}, \quad x\ge 0,\lambda> 0,
	\eeqlb
	where $t\mapsto v_t(\lambda)$ is the unique strictly positive solution to the differential equation
	\beqlb\label{CB-equation}
	\frac{\partial}{\partial t} v_t(\lambda)= -\itPsi(v_t(\lambda)), \quad v_0(\lambda)= \lambda.
	\eeqlb
	It is clear that $Q_t(0, \{0\})=1$ for all $t\ge 0$, which implies that $0$ is an absorbing state.
	The branching mechanism $\itPsi$ is said to be \textit{subcritical}, \textit{critical} or \textit{supercritical} according as $\itPsi'(0+)>0$, $\itPsi'(0+)=0$ or $\itPsi'(0+)<0$, respectively. In the sub-critical or critical situation, the CB-process will be surely absorbed by $0$ in finite time if and only if   Grey's condition is satisfied; see Condition \ref{con1}.

	The quasi-stationary distribution properties of CB-processes are very similar to those of Bienaym\'e--Galton--Waston branching processes. Lambert \cite{Lam07} pointed out that under Condition \ref{con1}, any subcritical CB-process has an infinite number of QSDs, and  characterized all of them. His study is based on the branching property, i.e., for any $t\geq 0$ and $x,y\ge0,$
	\beqlb\label{bp}
	Q_t(x, \cdot)*Q_t(y, \cdot)=Q_t(x+y, \cdot),
	\eeqlb
	where ``\,$*$'' denotes the convolution operator.
	This property means that different particles in population act independently of each other and allow the CB-processes to have rich mathematical structures.  However, due to limited resources or other mechanisms, the branching property \eqref{bp} is unlikely to be realistic in the real world.  Those considerations motivate the study of  CB-processes with competition.
	
	A well-known deterministic model that can be used to model competition is the so-called \textit{logistic growth model}, which is defined by the following ordinary differential equation:
	\beqnn
	\d Z_t=(bZ_t -a Z^2_t )\,\d t.
	\eeqnn
	Here the quadratic term can be interpreted ecologically, as it quantifies the negative interactions between each pair of individuals within the population. A stochastic analogue of this model was studied in Lambert \cite{Lam05} where  a \textit{logistic growth branching process} was introduced  combining the features of the logistic growth model and the CB-process.  The process was constructed in \cite{Lam05} by a random time change of a spectrally positive Ornstein--Uhlenbeck type process.
	Recently, Berestycki et al.\ \cite{BFF18} generalised the logistic growth branching process by considering a more general form of competition. Let $(\Omega, \mathscr{F}, \mathscr{F}_t, \Pp)$ be a filtered probability space satisfying the usual hypotheses. Suppose that $W(\d s,\d u)$ is a Gaussian white noise on $(0,\infty)^2$ based on $\d s\d u$, $M(\d s,\d z,\d u)$ is a Poisson random measure on $(0,\infty)^3$ with intensity $\d s\mu(\d z)\d u$, and  $\tilde{M}(\d s,\d z,\d u)$ is the compensated measure. Berestycki et al.\ \cite{BFF18} constructed the \textit{continuous-state branching process with competition} (CBC-process for short) as the unique strong solution to the following stochastic differential equation (SDE for short)
	\beqlb\label{sde1}
	Y_t \ar=\ar Y_0 + \sqrt{2c}\int_0^t\int^{Y_{s-}}_0 W(\d s, \d u) - \int_0^t [bY_{s}+g(Y_{s})]\d s
	\cr
	\ar\ar
	+ \int_0^t\int_0^1 \int_0^{Y_{s-}} z \tilde{M}(\d s,\d z,\d u) + \int_0^t\int_1^\infty \int_0^{Y_{s-}} z M(\d s,\d z,\d u),
	\eeqlb
	which generalized SDE construction for CB-process established by Dawson and Li \cite{DaL12}.
	Here, $b\in\mathbb{R}$ and $g$ is a continuous and non-decreasing function with $g(0)=0$, describing the intensity of competition between individuals. In particular, when $g(x)=ax^2$ for some $a>0$, the solution corresponds to the logistic branching process. When $g(x)\equiv0$, the solution reduces to a CB-process with the transition semigroup satisfying \eqref{CB-semigr}. This is the only situation where the branching property \eqref{bp} holds. If $g(x)$ is not identically $0$, the branching property and traditional methods  fail, and we need to find a new approach to study this process.

	Throughout this paper, we always assume that \eqref{sde1} has a unique non-explosive strong solution denoted by $(Y_t)_{t\ge0}$; see Theorem \ref{thm1} for related discussions. Recall that $\tau^-_0:=\inf\{t\ge0: Y_t=0\}$ is the extinction time
	 with the convention $\inf\emptyset=\infty$. Since all of the coefficients in equation \eqref{sde1} degenerate at state $0$, we know $0$ is  an absorbing state for $(Y_t)_{t\ge0}$; that is, $Y_t=0$ a.s. for all $t\ge \tau^-_0$, which agrees with the intuition  that whenever the population size attains $0$, nothing happens anymore, and the population size stays at $0$ forever. We are going to establish the uniqueness and existence of the QSD of $(Y_t)_{t\ge0}$, as well as the exponential convergence of conditional distribution.

	To the best of our knowledge, the existing studies on this topic only concern the case  that the  jump terms in \eqref{sde1} vanish. In this instance, the solution reduces to the so-called \textit{generalized Feller diffusion process.} Cattiaux et al. \cite{CCLMMS} established the  sufficient conditions  for the existence and uniqueness of the QSD for this process. In their proofs, spectral theory for self-adjoint operators plays an important role. However, this approach is not generally applicable since the Markov process corresponding to the solution to \eqref{sde1} is no longer symmetric when jump terms present.
More recently, the QSD problem of this diffusion process was revisited by Champagnat and Villemonais \cite{ChV15} through a probability-based approach established in \cite{ChV16}. Nevertheless, their method requires an accurate estimate for the hitting probability, which seems to be difficult to establish when jump terms exist.
	
	To present our main result, we need the following two conditions.
	\bgcondition\label{con1}
	{\rm(Grey's condition)}
	There exists a constant $\theta>0$ such that $\itPsi(\lambda)>0$ for all $\lambda>\theta$ and
	\beqlb\label{Grey}
	\int^\infty_\theta \frac{\d \lambda}{\itPsi(\lambda)}<\infty.
	\eeqlb
	\edcondition
	\begin{rem1}\rm \begin{itemize}
\item[{\rm (i)}]
		Condition \ref{con1} is a necessary and sufficient condition for a CB-process with branching mechanism $\itPsi$ to hit $0$ in finite time with strictly positive probability; see Grey \cite{Grey}.
\item[{\rm (ii)}]
		It is clear that Condition \ref{con1} is satisfied if $c>0$ in (\ref{bran-mech}).  On the other hand, if there exist constants $\beta\in(1, 2)$ and $C> 0$ such that
		\beqnn
		\mu(\d z)\ge \I_{\{0<z\le 1\}}C z^{-(1+\beta)} \d z,
		\eeqnn
		then Condition \ref{con1} is also satisfied. To see this, we have
		\beqnn
		\itPsi(\lambda)\ar=\ar b\lambda + c\lambda^2 + \int_0^\infty \big(\e^{-\lambda z}-1+\lambda z\I_{\{z\le 1\}}\big)\mu(\d z)\\
		\ar\ge\ar
		b\lambda +C\int_0^1 \big(\e^{-\lambda z}-1+\lambda z\big)z^{-(1+\beta)}\d z-\mu[1, \infty)\\
		\ar\ge\ar
		b\lambda +C\int_0^\infty \big(\e^{-\lambda z}-1+\lambda z\big)z^{-(1+\beta)}\d z-C\lambda \int_1^\infty z^{-\beta}\,\d z-\mu[1, \infty)\\
		\ar=\ar
		b\lambda +\frac{C\Gamma(2-\beta)}{\beta-1}\lambda^{\beta}-\frac{C\lambda}{\beta-1}-\mu[1, \infty).
		\eeqnn
		Hence, Condition $\ref{con1}$ is satisfied.\end{itemize}
	\end{rem1}
	
	Note that Condition \ref{con1} is necessary for CB-processes to have a QSD.
	This restriction unfortunately excludes some interesting branching mechanisms. An important example  is the so-called \textit{Neveu's branching process} (with branching mechanism $\itPsi(\lambda)=\lambda\log\lambda$ and L\'evy measure $\mu(\d z)=z^{-2} \d z$), which has rich mathematical structures and has found applications in various important area. In particular, Bertoin and Le Gall \cite{BeL00} provided a precise description of the connection between the so-called Neveu's CB-process and the coalescent processes introduced by Bolthausen and Sznitman \cite{BoS98} in their study of spin glasses. Berestycki et al. \cite{BBS13} proved that Neveu's CB-process may arise as a limit theorem of certain rescaled particle systems.
	
	However, Condition \ref{con1} is not necessary for CBC-process to have a QSD.  For the CBC-process $(Y_t)_{t\ge0}$ with Neveu's type branching mechanism to have a unique QSD, we need an alternative condition on the competition mechanism near $0$. Since  Condition \ref{con1} holds for $c\neq 0$, we only need to impose a condition for $c=0$.
	\bgcondition\label{con2}
 $c=0$ and  there exist constants $\theta\in (0,1)$ and $C>0$ such that $$\liminf_{x\to 0}g(x)x^{-\theta}>0\qquad\mbox{and}\qquad \mu(\d z)\ge \I_{\{0<z\le 1\}}C z^{-2} \d z.$$
	\edcondition

	We next present the main result in this paper.
	
	\bgtheorem\label{QSD} Suppose that one of Conditions~$\ref{con1}$ and $\ref{con2}$ holds and in addition  there are constants $c_0>0$ and $
	\alpha\in (0, 2)$ such that $$\frac{\I_{\{z>1\}}\mu(\d z)}{\d z}\le c_0  z^{-1-\alpha}\d z.$$ Recall that function $g$ is  continuous and non-decreasing with $g(0)=0$. For any increasing and strictly positive function $\varphi\in C^1[0,\infty)$ satisfying $\int_1^\infty \frac{1}{r\varphi(r)}\,\d r<\infty$, if the following condition further holds
	$$\begin{cases}
		\displaystyle\lim_{x\to \infty} \frac{g(x)}{x\varphi(x)}=\infty,&\quad \alpha\in (1,2),\\
		\displaystyle\lim_{x\to \infty} \frac{g(x)}{x\varphi(x)\log x}=\infty,&\quad \alpha=1,\\
		\displaystyle\lim_{x\to \infty} \frac{g(x)}{x^{2-\alpha}}=\infty,&\quad \alpha\in(0,1),\end{cases}$$
		then the solution $(Y_t)_{t\geq 0}$ to \eqref{sde1} has a unique quasi-stationary distribution {$\pi$}. Moreover,
there exists a constant $\lambda>0$
 such that for any Borel set $A\subset (0, \infty)$,
$t>0$ and any initial distribution $\nu$ on $(0, \infty)$,
	\[|\mbb{P}_{\nu}(Y_t\in A|\tau^-_0>t)-\pi(A)|\le C(\nu)\e^{-\lambda t},\] where $\nu\mapsto C(\nu)$ is a non-negative function.
	\edtheorem
	
\subsection{Strategy of the approach}
	The proof of Theorem \ref{QSD} relies on the recent paper by Guillin et al. \cite{GNW}, which establishes general results regarding the existence and uniqueness of the QSD of strongly Feller Markov processes.
		Denote by $(P_t)_{t\ge0}$ the transition semigroup of $(Y_t)_{t\ge 0}$. Let $b\mcr{B}[0,\infty)$ be the space of all bounded Borel functions on $[0, \infty)$, and let $C_b[0,\infty)$ be the space of continuous  bounded functions on $[0, \infty)$.  For $T>0$, denote by $D([0, T], \mbb{R}_+)$ the space of $\mbb{R}_+$-valued c\'adl\'ag paths defined on $[0, T]$ and equipped with the Skorokhod topology. We next introduce five conditions to be verified under Condition \ref{con1} or \ref{con2}.
	
	\noindent{\bf (C1) (Strong Feller property)}  For any $t>0$ and $f\in b\mcr{B}[0,\infty)$,
	$x\mapsto P_tf(x)$ is continuous on $[0, \infty).$
	
	\noindent{\bf (C2) (Trajectory Feller property)}
	For any $T>0$, $x\mapsto \mbb{P}_x((Y_t)_{t\in[0, T]}\in \cdot)$ is continuous in the sense of weak convergence.
	
	From It\^{o}'s formula and \eqref{sde1} one can check that for any $f\in C^2_b[0, \infty)$ and $x\in [0,\infty)$,
	\[\mbb{E}_xf(Y_t)=f(x)+\mbb{E}_x\left[\int^t_0 Lf(Y_s) \d s\right],\]
	where
	\beqlb\label{gen-CBC}
	Lf(x) \ar:=\ar  -[bx+g(x)]f^\prime(x)+cxf^{\prime\prime}(x)\cr
	\ar\ar
	+\,  x\int_0^\infty \big[f(x+z) - f(x) - zf^\prime(x)\I_{\{z\le 1\}}\big] \mu(\d z).
	\eeqlb

	\noindent{\bf (C3) (Lyapunov condition)} There exist a function $W\in C^2_b[0, \infty)$
	such that $W(x)\ge 1$ for all $x\in [0, \infty)$, and there are two sequences of positive constants $(r_n)_{n\ge1}$ and $(b_n)_{n\ge1}$, where $r_n\to\infty$ as $n\to \infty$, and an increasing sequence of compact subsets $(K_n)_{n\ge1}$ of $[0, \infty),$ such that for all $n\ge1$,
	\[-LW(x)\ge r_n W(x)-b_n\I_{K_n}(x),
	\quad x\in [0,\infty).\]
	
	Let $(P^0_t)_{t\ge 0}$ be the transition semigroup of the killed process $(Y_t)_{0\le t<\tau^-_0}$; that is, for any $t\ge0$ and $f\in b\mcr{B}[0,\infty)$,
	\[P_t^0f(x):=\mbb{E}_x[\I_{\{t< \tau^-_0\}}f(Y_t)],\quad x\in (0,\infty).\] Denote by $P_t^0(x, \d y)$ the transition kernel of the killed process $(Y_t)_{0\le t<\tau^-_0}$.
	
	\noindent{\bf (C4) (Weak Feller property)} For a measure-separable class $\mcr{A}$ of continuous bounded functions $f$ with support contained in $(0, \infty)$, $x\mapsto P^0_t f(x)$ is continuous on $(0, \infty).$
	
	\noindent{\bf (C5) (Irreducibility)} For all $x, t>0$ and non-empty open set $O\subset (0, \infty)$,
	$P_t^0(x, O)>0.$ In addition, there exists $x_0>0$ such that $\mbb{P}_{x_0}(\tau^-_0<\infty)>0.$
	
	It follows from \cite[Theorem 2.2]{GNW} that, if all of the conditions {\bf (C1)}--{\bf (C5)} hold, then the result of Theorem \ref{QSD} is valid. Therefore, before showing Theorem \ref{QSD}, we first present results concerning  (C1) to (C5), respectively.

	\bgproposition\label{pc1} Suppose that either Condition~$\ref{con1}$ is satisfied or there exist constants $\theta\in (0,1)$, $\beta\in(0, 2)$ and $C> 0$ such that $$\liminf_{x\to 0}g(x)x^{-\theta}>0\qquad\mbox{and}\qquad \mu(\d z)\ge \I_{\{0<z\le 1\}}C z^{-(1+\beta)} \d z.$$
	Then $(P_t)_{t\ge 0}$ has the strong Feller property {\bf (C1)}.
	\edproposition
	
	\bgproposition\label{pc2}
	The process $(Y_t)_{t\ge 0}$ has the trajectory Feller property {\bf (C2)}.
	\edproposition
	
	\bgproposition\label{pc3}
	Assume that $$\frac{\I_{\{z>1\}}\mu(\d z)}{\d z}\le c_0  z^{-1-\alpha}\d z$$ for some  $\alpha\in (0,2)$ and $c_0>0$. Let $\varphi\in C^1[0,\infty)$ be  an increasing and strictly positive function satisfying $\int_1^\infty \frac{1}{r\varphi(r)}\,\d r<\infty.$ If any of  the following conditions further holds
	$$\begin{cases}
		\displaystyle\lim_{x\to \infty} \frac{g(x)}{x\varphi(x)}=\infty &\quad\text{if}\quad \alpha\in (1,2),\\
		\displaystyle\lim_{x\to \infty} \frac{g(x)}{x\varphi(x)\log x}=\infty &\quad\text{if}\quad \alpha=1,\\
		\displaystyle\lim_{x\to \infty} \frac{g(x)}{x^{2-\alpha}}=\infty &\quad\text{if}\quad \alpha\in(0,1),\end{cases}$$ then the Lyapunov condition {\bf (C3)} is satisfied.
	\edproposition
	
	\bgproposition\label{pc4}
	Suppose that $(P_t)_{t\ge 0}$ has the strong Feller property. Then the associated transition semigroup $(P^0_t)_{t\ge 0}$  has the weak Feller property {\bf (C4)}.
	\edproposition
	
	\bgcondition\label{nontri}
	Either $c>0$ or $\int^1_0 z \mu(\d z)=\infty$.
	\edcondition
	Indeed, Condition \ref{con1} implies  Condition \ref{nontri}. If $c=0$ and $\int^1_0 z \mu(\d z)<\infty$, then we have
	\beqnn
	\itPsi(\lambda)
	\ar=\ar b\lambda + \int_0^\infty \big(\e^{-\lambda z}-1+\lambda z\I_{\{z\le 1\}}\big)\mu(\d z)\\
	\ar\le\ar
	b\lambda +\int_0^1 \big(\e^{-\lambda z}-1+\lambda z\big)\mu(\d z)\\
	\ar\le\ar
	b\lambda +\lambda\int_0^1 z\mu(\d z),
	\eeqnn
	and hence Condition \ref{con1} is not satisfied.
	
	\bgproposition\label{pc5}
	Suppose that Condition $\ref{nontri}$ holds.
	Then the irreducibility condition {\bf (C5)} is satisfied.
		\edproposition

With all the propositions above at hand, we can give the

\medskip\noindent\textit{Proof of Theorem~$\ref{QSD}$.~}
Applying Propositions
\ref{pc1}, \ref{pc2}, \ref{pc3}, \ref{pc4} and \ref{pc5}, we know that all the conditions {\bf (C1)}--{\bf (C5)} hold under Condition \ref{con1} or Condition \ref{con2}.
Then, it follows from \cite[Theorem 2.2]{GNW} that the statements in Theorem \ref{QSD} are valid.\qed

One of the main difficulties  comes from the fact that the coefficients tend to $0$ when the process is close to $0$, which is clear from \eqref{sde1}. We shall mention two novel points in the proof of the propositions above. One is that, in order
to show the strong Feller property for the CBC-process, we adopt a localization argument that reduces the required condition; see Subsection \ref{section2.2}. The other is that, in showing the irreducibility, we also apply the Lamperti transform and comparison principle  to reduce the problem to
that on L\'evy process; see Subsection  \ref{section2.4}.

	\section{Proof}
	
	In this section, we are going to prove Propositions \ref{pc1}--\ref{pc5}. For the sake of completeness, we begin with the non-explosion of the unique strong solution to the SDE \eqref{sde1}.
	
	\subsection{Uniqueness of the strong solution}
Let $L$ be the operator defined by \eqref{gen-CBC}.
	Let $\mcr{D}(L)$ denote the linear space consisting of twice continuously differentiable functions on $[0, \infty)$, such that the integral on the right-hand side of \eqref{gen-CBC} converges, and $Lf$ is a continuous function on $[0, \infty)$.
	
	\bgtheorem\label{thm1} Suppose that there exist constants $C_1, C_2>0$ and a non-negative function $V\in\mcr{D}(L)$ such that $\lim_{x\to\infty}V(x)=\infty$ and
	\beqlb\label{Nonexp-cond}
	LV(x)\le C_1+C_2V(x),\qquad
	x\in [0,\infty).
	\eeqlb
	Then, for any $y\ge 0$, there is a pathwise unique non-explosive strong solution $\{Y_t: t\ge 0\}$ to \eqref{sde1} with $Y_0=y$.
	\edtheorem
	\proof
	By Theorem~2.5 of Dawson and Li \cite{DaL12}, for each $n\ge 1$ there is a pathwise unique strong solution $\{y^{(n)}_t: t\ge 0\}$ to
	\beqlb\label{sde-CBIC_n}
	y_t \ar=\ar y + \int_0^t\int_0^{{\red y_{s}}}W(\d s,\d u) + \int_0^t [-by_s-g(y_s)] \d s \cr
	\ar\ar
	+ \int_0^t\int_0^1 \int_0^{y_{s-}} z \tilde{M}(\d s,\d z,\d u) + \int_0^t\int_1^\infty \int_0^{y_{s-}} (z\land n) M(\d s,\d z,\d u).
	\eeqlb
	Let $\zeta_n:= \inf\{t\ge 0: y^{(n)}_{t}\ge n\}$. From the pathwise uniqueness for \eqref{sde-CBIC_n} it follows that $\zeta_n$ is nondecreasing in $n\ge 1$ and $y^{(n+1)}_t= y^{(n)}_t$ for $0\le t< \zeta_n$. Clearly, the pathwise unique solution $\{y_t: t\ge 0\}$ to \eqref{sde1} is given by $y_t= y^{(n)}_t$ for $0\le t< \zeta_n$ and $y_t= \infty$ for $t\ge \zeta:= \lim_{n\to \infty}\zeta_n$.
	
	We next prove that the solution $\{y_t: t\ge 0\}$ is non-explosive, i.e., the explosion time $\zeta$ is infinite almost surely.
	By It\^{o}'s formula, we have \beqlb\label{LV}
	V(y_{t\land \zeta_n})= V(y) + \int_0^{t\land \zeta_n} LV(y_s) \d s + M_{n}(t),
	\eeqlb
	where
	\beqlb\label{M_n(t)}
	M_{n}(t)\ar=\ar \int_0^{t\land \zeta_n}\int_0^{y_{s}}V'(y_{s-})W(\d s,\d u)\cr
	\ar\ar
	+ \int_0^{t\land \zeta_n}\int_0^\infty\int_0^{y_{s-}} (V(y_{s-}+z)- V(y_{s-}))\tilde{M}(\d s,\d z,\d u)
	\eeqlb
	is a martingale.
	Taking  expectations  on both sides of \eqref{LV} and applying \eqref{Nonexp-cond}, we obtain
	\beqnn
	\mbb{E}[V(y_{t\land \zeta_n})]\le V(y) + \mbb{E}\left[\int^t_0 \left(C_1+C_2V(y_{s\land \zeta_n}\right)\d s\right].
	\eeqnn
	By Gronwall's lemma, it follows that, for any $T>0$ and $t\in (0,T]$,
	\beqnn
	\mbb{E}[V(y_{t\land \zeta_n})]\le (V(y)+C_1T) \e^{C_2t}.
	\eeqnn
	Letting $n\to\infty$ in the above and noting that $V(x)\to\infty$ as $x\to\infty$, we know that $\mbb{P}(\zeta\le t)= 0$.  By the arbitrariness of $t$ and $T$, we have $\zeta= \infty$ almost surely.\qed

	\subsection{Strong Feller property}\label{section2.2}
	
	The proof of the strong Feller property  under Condition \ref{con1} or Condition  \ref{con2} is mainly based on the probabilistic coupling method. We shall first prove the desired result  under Condition \ref{con2}, and then deal with Condition \ref{con1}.  In Li et al. \cite{LLWZ22}, we have constructed a two-dimensional \textit{Markov coupling process} $\{(X_t,Y_t): t\ge 0\}$ on $D:= \{(x,y): x\ge y\ge 0\}$, where both  $(X_t)_{t\ge 0}$ and $(Y_t)_{t\ge 0}$  are Markov processes with the same transition semigroup $(P_t)_{t\ge 0}$. The  coupling process satisfies $X_{T+t}= Y_{T+t}$ for every $t\ge 0$, where
	$$T:= \inf\{t\ge 0: X_t= Y_t\}$$ 	is the \textit{coupling time}.
	We next give a characterization for the generator of the Markov coupling process $\{(X_t,Y_t): t\ge 0\}$.
	Let $\itDelta= \{(z,z): z\ge 0\}\subset D$ and $\itDelta^c= D\setminus \itDelta$.
	Let $\mu_a= [\mu\land(\delta_a*\mu)]/2$ for any $a\in\R$. Given a function $F$ on $D$ which is twice continuously differentiable on $\itDelta^c$, we write for any $(x,y)\in \itDelta^c$,
	\beqlb\label{CoupleGen}
	\tilde{L}F(x,y)\ar:=\ar \big[-bx-g(x)\big]F'_x(x,y) + \big[-by-g(y)\big]F'_y(x,y) \ccr
	\ar\ar
	+\, cxF''_{xx}(x,y)+cyF''_{yy}(x,y) - 2cyF''_{xy}(x,y)\cr
	\ar\ar
	+\, (x-y)\int_0^\infty\big[F(x+z,y)-F(x,y)-F_x'(x,y)z\I_{\{z\le 1\}}\big] \mu(\d z) \cr
	\ar\ar
	+\,y\int_0^\infty \big[F(x+z,y+z)-F(x,y)-(F_x'+F_y')(x,y)z\I_{\{z\le 1\}}\big] \mu(\d z)\ccr
	\ar\ar+\, y\int_0^\infty \big[F(x+z,2y+z-x)-F(x+z,y+z)\big] \mu_{x-y}(\d z) \ccr
	\ar\ar
	+\,y\int_0^\infty \big[F(x+z,x+z)-F(x+z,y+z)\big] \mu_{y-x}(\d z).
	\eeqlb
	Let $\mcr{D}(\tilde{L})$ denote the linear space consisting of the functions $F$ such that the integrals in \eqref{CoupleGen} are finite and $\tilde{L}F$ is locally bounded on compact subsets of $\itDelta^c$. We call $(\tilde{L},\mcr{D}(\tilde{L}))$ the \textit{coupling generator} of the CBC-process $\{Y_t:t\ge0\}$.
	Let
	$$\zeta_n^*:= \inf\{t\ge 0: X_t\ge n\,\,\text{ or}\,\, X_t-Y_t\le 1/n\}$$
	and assume that $(X_0, Y_0)=(x, y)\in \itDelta^c$. Then, it follows from Li et al. \cite[Theorem 4.1]{LLWZ22} that for any $n\ge 1$ and $F\in \mcr{D}(\tilde{L})$
	\beqlb\label{Ito-f(X_t,Y_t)}
	F(X_{t\land \zeta_n^*},Y_{t\land \zeta_n^*}) = F(x,y) + \int_0^{t\land \zeta_n^*} \tilde{L}F(X_s,Y_s) \d s + M_n(t),
	\eeqlb
	where $\{M_n(t): t\ge 0\}$ is a martingale.

	Note that when $F(x, y)=f(x-y)\I_{\{(x, y)\in \itDelta^c\}}$ for some
	$f\in \mcr{D}(L)$, it is clear that $F\in \mcr{D}(\tilde{L}).$ In this case, we also write $\tilde{L}f(x,y):=\tilde{L}F(x, y)$ for convenience. Moreover, for any $(x,y)\in \itDelta^c$,
	\beqlb\label{CLf}
	\tilde{L}f(x,y)\ar=\ar -\big[b(x-y)+g(x)-g(y)\big]f'(x-y)\ccr
	\ar\ar
	+\, c(x-y)f''(x-y)+4cyf''(x-y)\cr
	\ar\ar
	+\, (x-y)\int_0^\infty\big[f(x+z-y)-f(x-y)-f'(x-y)z\I_{\{z\le 1\}}\big] \mu(\d z) \ccr
	\ar\ar+y[f(2x-2y)-2f(x-y)] \mu_{x-y}(\mbb{R}_+)\ccr
	\ar=\ar
	-[g(x)-g(y)]f'(x-y)+4cyf''(x-y)\ccr
	\ar\ar+\,y[f(2x-2y)-2f(x-y)] \mu_{x-y}(\mbb{R}_+)+Lf(x-y),
	\eeqlb
	where we used the fact that  $\mu_{x-y}(0,\infty)= \mu_{y-x}(0,\infty)$.
	
	Before proving the main result in this part, we first establish the following two useful lemmas.

	\bglemma\label{Lemma-uniform}
	Given two intervals $(A, B)$ and $(A', B')$ such that $0<A<A'<B'<B<\infty$, let $(Y_t)_{t\ge0}$  be a CBC-process with initial value $y\in (A', B')$, and define  $S:=\inf\{t>0: Y_t\notin (A, B)\}$. Then, there is a constant $C>0$ so that for any $t\in (0,1]$ and any $y\in (A', B')$,
	\[
	\mbb{P}_y({S} 
<t)\le Ct^{1/2}.
	\]
	\edlemma
	\proof
	Recall that  $(Y_t)_{t\ge0}$ is the strong solution to the equation \eqref{sde1}. For $t>0$, define
	$$ M_{t}
	=\sqrt{2c}\int_0^{t}\int^{Y_{s}}_0 W(\d s, \d u)+\int_0^{t}\int_0^1\int_0^{{Y}_{s-}} z \tilde{M}(\d s,\d z,\d u).$$
	Then,
	\begin{equation*}
		M_{t\wedge S}=Y_{t\wedge S}-Y_0+Y^{(1)}_{t\wedge S}+Y^{(2)}_{t\wedge S},\quad t\ge0,
	\end{equation*}
	where
	\[
	Y^{(1)}_t=\int_0^{t} [bY_{s }+g(Y_{s })]\d s
	\]
	and
	\[
	Y^{(2)}_t=-\int_0^{t}\int_1^\infty\int_0^{Y_{s-}} z {M}(\d s,\d z,\d u).
	\]

	Let $l=(A'-A)\wedge(B-B')$. Then, we have
	\begin{equation}\label{ineq3}
		\begin{split}
			\mbb{P}_y(S<t)
			&\leq \mbb{P}_y(|Y_{t\wedge S}-y|\geq l)  \\
			&\leq\mbb{P}_y(|M_{t\wedge S}|>l/2 )+\mbb{P}_y(|Y^{(1)}_{t\wedge S}|>l/2)+\mbb{P}_y(Y^{(2)}_{t\wedge S}\neq 0).
		\end{split}
	\end{equation}
	
	It is clear that $(M^*_{t})_{t\ge0}:=(M_{t\wedge S})_{t\ge0}$ is a martingale with  quadratic variation process
	\[ \langle M^*\rangle_{t}=\int_0^{t\wedge S}2c{Y_s}\d s+\int_0^{t\wedge S}Y_s \d s\int_0^1 z^2\mu(\d z),\quad t\ge0.\]
	Using the Jessen inequality and the Markov inequality we have
	\begin{equation}\label{ineq4}
		\begin{split}
			\mbb{P}_y(|M_{t\wedge S}| >l/2)
			&\leq \frac{2}{l}\mbb{E}_y|M^*_{t}| \leq \frac{2}{l}(\mbb{E}_y|{M^*_{t}}|^2)^{1/2}=\frac{2}{l}(\mbb{E}_y\langle M^*\rangle_{t})^{1/2} \\
			&=\frac{2}{l}\left\{\mbb{E}_y\left[\int_0^{t\wedge S}2c {Y}_s \d s\right]+\mbb{E}_y\left[\int_0^{t\wedge S}{Y}_s\d s\int_0^1 z^2\mu(\d z)\right]\right\}^{1/2}\\
			&\le \frac{2}{l}\left(2c+\int_0^1 z^2\mu(\d z)\right)^{1/2}t^{1/2}B^{1/2}.
		\end{split}
	\end{equation}
	
	Since $x\mapsto g(x)$ is non-decreasing, we have
	\[
	|Y^{(1)}_{t\wedge S}|\le t(|b|B+g(B))
	\]
	and hence
	\beqlb\label{ineq5}
	\mbb{P}_y(|Y^{(1)}_{t\wedge S}|>l/2)\le 2t(|b|B+g(B))/l.
	\eeqlb

	Noting that $(Y^{(2)}_t)_{t\ge 0}$ is a pure jump process, and the waiting time before its first jump can be controlled by an exponential random variable, we have
	\beqlb\label{ineq6}
	\mbb{P}_y(Y^{(2)}_{t\wedge S}\neq 0)\le 1-\exp\left(-B\mu(1, \infty)t\right)\le B\mu(1, \infty)t.
	\eeqlb
	
	Finally, substituting \eqref{ineq4}, \eqref{ineq5} and \eqref{ineq6} into \eqref{ineq3} gives the desired result.
	\qed

	The following lemma provides a way to establish the
	strong Feller property by verifying a local inequality of the coupling generator $\tilde L$. We refer to Luo and Wang \cite[Theorem 5.1]{LuW19} for a similar result in which a universal inequality is required.
	\bglemma\label{l1}
	Given any interval $(A, B)$ with $0<A,B<\infty$, if there exist  constants $l, C>0$ and a non-negative  continuous function $\phi\in \mcr{D}(L)$  such that $\phi(0)=0$ and
	\beqlb\label{LY1}
	\tilde{L}{\itPhi(x, y)}\le -C, \qquad x, y\in (A, B) \hbox{ with }   l>x-y>0,
	\eeqlb
	{where $\itPhi(x, y):=\phi(x-y)\I_{\{(x,y)\in\itDelta^c\}}$}, then for any $t>0$ and  any $f\in b\mcr{B}[0,\infty)$, $x\mapsto P_tf(x)$ is uniformly continuous on $(A, B)$.
	\edlemma
	\proof
	Set
	\[T_n:=\inf\{t\ge0: X_t-Y_t\le 1/n\}\]
	and
	\[{S}^*:=\inf\{t>0: \text{either}\,\, X_t\notin (A, B)~\mbox{or}~Y_t\notin (A, B)\}.\]
	{Then the coupling time $T=\lim_{n\to\infty} T_n.$} Let $(A', B')$ be a sub-interval satisfying $A<A'<B'<B$.
	Assume that $X_0=x$ and $Y_0=y$ with $A'<y<x<B'$. Take $n$ large enough such that $x-y>1/n$. Then, by \eqref{Ito-f(X_t,Y_t)}, we have
	\beqnn
	\phi(X_{t\wedge T_{n}\wedge { S}^*}-Y_{t\wedge T_{n}\wedge {S}^*})={ \phi(x-y)}+\int^{t\wedge T_{n}\wedge { S}^*}_0\tilde{L}\phi(X_s, Y_s)\d s+M_{t\wedge T_{n}\wedge {\red S}^*}.
	\eeqnn
	
	Since $\tilde{L}\phi(x,y)$ is bounded for $x- y\in [1/n, l]$, we see that $(M_{t\wedge T_{n}\wedge S^*})_{t\ge 0}$ is a martingale.
	By taking expectations, and using the fact that
	$X_t>Y_t$ for $t\leq T_n$, and (\ref{LY1}), we see that
	\beqnn
	0\ar\le\ar\mbb{E}[\phi(X_{t\wedge T_{n}\wedge S^*}-Y_{t\wedge T_{n}\wedge S^*})]\\
	\ar=\ar
	\phi(x-y)+\mbb{E}\left[\int^{t\wedge T_{n}\wedge S^*}_0\tilde{L}\phi(X_s, Y_s)\d s\right]\\
	\ar\le\ar
	\phi(x-y)-C \mbb{E}[t\wedge T_{n}\wedge S^*].
	\eeqnn
	Letting $t\to\infty$ and $n\to\infty$, we obtain
	\beqlb\label{LY2}
	\mbb{E}[T\wedge S^*]\le\frac{\phi(x-y)}{C}.
	\eeqlb
	It follows that
	\[\mbb{P} (T\wedge S^*>t)\leq \frac{\phi(x-y)}{Ct}\]
	and
	\[\phi(x-y)^{1/2}\mbb{P} (\phi(x-y)^{1/2}<S^*< T)\leq \frac{\phi(x-y)}{C}.\]

	Combining  all the estimates above, we obtain
	\begin{equation}\label{1}
		\begin{split}
			{\mbb{P}}(T>t)
			&\le \mbb{P}(T\wedge S^*>t)+\mbb{P}(S^*< T) \\
			&\le\mbb{P}\{T\wedge S^*>t\}+\mbb{P} (\phi(x-y)^{1/2}<S^*< T)\\
			&\quad +\mbb{P}(S^*\leq \phi(x-y)^{1/2})\\
			&\le \frac{\phi(x-y)}{tC}+\frac{\phi(x-y)^{1/2}}{C}+\mbb{P}(S^*\leq \phi(x-y)^{1/2}).
		\end{split}
	\end{equation}
Define $S^X=\inf\{t>0:X_t\notin (A,B)\}$ and $S^Y=\inf\{t>0:Y_t\notin (A,B)\}$.
Applying Lemma \ref{Lemma-uniform}, \eqref{1} and the fact that
	\[\mbb{P}(S^*\leq \phi(x-y)^{1/2})
	\le \mbb{P}(S^X\leq \phi(x-y)^{1/2})
	+\mbb{P}(S^Y\leq \phi(x-y)^{1/2}),
	\]
	we obtain $$\lim_{x-y\rightarrow 0+}\mbb{P}(S^*\leq  \phi(x-y)^{1/2})=0$$ uniformly for $A'<y<x<B'.$

	Thus, by the arbitrariness of $(A', B')$, we know for any $f\in b\mcr{B}[0,\infty)$, $t>0$ and $A<y<x<B$,
	\beqnn
	|P_tf(x)-P_tf(y)|\ar=\ar
	|\mbb{E}[(f(X_t)-f(Y_t))\I_{\{T>t\}}]|\\
	\ar\le\ar 2\|f\|_\infty\mbb{P}(T>t)\to 0\,\,\,\text{as}\,\,\, x-y\to 0+.
	\eeqnn Hence, the desired assertion follows. \qed

	\bgcondition\label{con2'} {\rm(Fluctuation condition)}
	There exist constants $\beta\in(0, 2)$ and $C> 0$ such that
	\beqlb\label{mux-cond}
	\mu(\d z)\ge \I_{\{0<z\le 1\}}C z^{-(1+\beta)} \d z.
	\eeqlb
	\edcondition

\bgproposition\label{p1}
Suppose that Condition~$\ref{con2'}$ is satisfied. Then, for any $f\in b\mcr{B}[0,\infty)$ and any $t>0$, the mapping $x\mapsto P_tf(x)$ is continuous on $(0, \infty)$.
\edproposition
\proof
It follows from \cite[Example 1.2]{LuW19} that
there are some $C_*, \kappa>0$ so that for all $x\in(0, \kappa)$
\beqlb\label{est24}\mu_x(\mbb{R_+})> C_*x^{-\beta}.
\eeqlb
We take $\rho=(\beta\wedge1)/2$.
Define
\beqlb\label{phi}
\phi(r):= 1-\e^{-r^\rho}, \quad r\ge 0.
\eeqlb

It is clear that for any $r>0$, \[\phi(0)=0,\quad\phi'(r)=\rho r^{\rho-1} \e^{-r^\rho}>0, \quad
\phi''(r)=-\rho r^{\rho-2} \e^{-r ^\rho}[1-\rho+\rho r^\rho]<0\]
and
$$\phi^{(3)}(x)=\rho(2-\rho)r^{\rho-3}\e^{-r ^\rho}[1-\rho+\rho r^\rho]+\rho^2r^{2\rho-3}\e^{-r ^\rho}[1-2\rho+\rho r^\rho]\ge0.$$  In particular, for all $x\in (0,1/2)$,
\begin{align*}&\int_0^\infty\big[\phi(x+z)-\phi(x)-\phi'(x)z\I_{\{z\le 1\}}\big] \mu(\d z)\\
	&\qquad
\le\int_0^x\big[\phi(x+z)-\phi(x)-\phi'(x)z\big] \mu(\d z)+ \mu(1, \infty)\\
	&\qquad\le \frac{1}{2}\phi''(2x)\int_0^xz^2\mu(\d z)+ \mu(1, \infty),
\end{align*} where in the last inequality we used the mean value theorem and the fact that $\phi^{(3)}\ge0$.

Recalling the definition of $L$ in \eqref{gen-CBC}, we have for all $x\in (0,1/2)$,
\begin{equation}\label{Lphi}\begin{split}
		L\phi(x)
		&\le
		(-bx-g(x))\phi'(x)+\left[c+\frac{1}{2}\int_0^xz^2\mu(\d z)\right]x\phi''(2x)+ x \mu(1, \infty)\\
		&\le |b|\rho x^\rho \e^{-x^\rho} +x \mu(1, \infty).
\end{split}\end{equation}

On the other hand, by the mean value theorem again and the facts that $\phi(0)=0$, $\phi''\le 0$ and $\phi^{(3)}\ge0$,
we know
\begin{align*}\phi(2x)-2\phi(x)=&\int_0^{2x}\phi'(u)\,\d u-2\int_0^x\phi'(u)\,\d u=\int_0^x\phi'(x+u)\,\d u-\int_0^x\phi'(u)\,\d u\\
	\le& \int_0^x \phi'(2x)x\,\d u=\phi''(2x)x^2.\end{align*}

Fix $0<A<B$ and $0<l<1/2$. Apply \eqref{CLf} with $f$ replaced by $\phi$. Since the function $g$ is non-decreasing, for any $x,y\in (A, B)$ with $0<x-y<l$, we have
\beqlb\label{CLphi}
\tilde{L}\phi(x, y)\ar\le\ar
|b|\rho (x-y)^\rho \e^{-(x-y)^\rho}+(x-y) \mu(1, \infty)\ccr   \ar\ar-4cy\rho(1-\rho)(x-y)^{\rho-2}\e^{-(x-y)^\rho}\ccr
\ar\ar-2y\rho (x-y)^2\phi''(2(x-y))\mu_{x-y}(\mbb{R_+})\ccr
\ar\le\ar
|b|\rho l^\rho+l \mu(1, \infty)\ccr
\ar\ar-2A\rho (1-\rho)2^{\rho-2}(x-y)^{\rho}\e^{- (2l)^\rho}\mu_{x-y}(\mbb{R_+}).
\eeqlb

It follows from \eqref{est24} that the last item  on the right-hand side of \eqref{CLphi} tends to $-\infty$ as $l\to0$.
Then, we can choose small enough $l>0$ such that $\tilde{L}\phi(x, y)<-1$ for all  $x,y\in (A, B)$ with $0<x-y<l$.

This along with Lemma \ref{l1} and the arbitrariness of the interval $(A, B)$ immediately gives the desired assertion.
\qed

\bgproposition\label{p2} Suppose that
$$\liminf_{x\to 0+}g(x)x^{-\theta}>0$$ for some $\theta\in (0,1)$. Then, for any $t>0$ and any function $f\in b\mcr{B}[0,\infty)$,
\[
\lim_{x\to0+} P_tf(x)=P_tf(0).
\] \edproposition

\begin{proof}
	Let $\phi$ be the function defined by \eqref{phi} with $\rho= (1-\theta)/2$.
	According to the first inequality in \eqref{Lphi}, we see that there is a sufficiently small $\delta\in (0, l)$ such that for any  $x\in(0,\delta)$,
	\[L\phi(x)<-1.\] Denote by $\{Y_t\}_{t\ge0}$ the solution to \eqref{sde1} with initial value $y>0$, and define
	$$\tau^+_{\delta}:=\inf\{t>0: Y_t>\delta\}\,\,\text{ and }\,\, \tau^-_\varepsilon:=\inf\{t>0: Y_t\leq \varepsilon\}\,\,\text{ for }\,\, \delta>\varepsilon\ge0.$$
	Similar to the proof of Theorem \ref{thm1}, applying It\^{o}'s formula we obtain
	\beqnn
	\phi(Y_{t\wedge\tau^-_\varepsilon\wedge\tau^+_{\delta}})=\phi(y)+\int^{t\wedge\tau^-_\varepsilon\wedge\tau^+_{\delta}}_0L\phi(Y_s)\d s+M_t,
	\eeqnn
	where $\{M_t:t\ge0\}$ is a martingale. Taking expectations on both sides of the above, we have
	\beqlb\label{mp1}
	0\le\mbb{E}_y[\phi(Y_{t\wedge\tau^-_\varepsilon\wedge\tau^+_{\delta}})]\le\phi(y)-\mbb{E}_y[t\wedge\tau^-_\varepsilon\wedge\tau_{\delta}]
	\eeqlb
	Letting $t\to\infty$ in the inequality above, we obtain
	\[
	\mbb{E}_y[\tau^-_\varepsilon\wedge\tau^+_{\delta}]\le\phi(y).
	\]
	Since $\phi$ is continuous and $\phi(0)=0$, using the Markov inequality, for any $t>0$, we have
	$$
	\mbb{P}_y \left(\tau^-_0\wedge\tau^+_{\delta}>t\right)\le \phi(y)/t\to0\quad\mbox{as}\quad y\to 0,
	$$
	On the other hand, using (\ref{mp1}) again gives
	$$
	\mbb{P}_y\left(\tau^+_{\delta}\le t\wedge \tau^-_0\right)\le\phi(y)/\phi(\delta)\to0\quad\mbox{as}\quad y\to 0.
	$$

	Putting two inequalities above together gives
	\beqlb\label{neweq1}
	\mbb{P}_y\left(\tau^-_0>t\right)\le\mbb{P}_y \left(\tau^-_0\wedge \tau^+_\delta>t\right)+\mbb{P}_y \left(\tau^+_\delta\le t\wedge\tau^-_0\right)\le \phi(y)/t+\phi(y)/\phi(\delta),
	\eeqlb
	which tends to $0$ as $y\to0$.
	Thus, for any $f\in b\mcr{B}[0,\infty)$ and $t>0$,
	\beqnn
	|P_tf(y)-P_tf(0)|\ar\le\ar
	|P_tf(y)-P_tf(0)|\\
	\ar\le\ar
	\left|\mbb{E}_y\left[\left(f(Y_t)-f(0)\right)\I_{\{\tau^-_0>t\}}\right]\right|\\
	\ar\le\ar 2\|f\|_\infty\mbb{P}_y(\tau^-_0>t)\to 0 \quad\text{ as}\quad y\to 0.
	\eeqnn The proof is complete. \qed\end{proof}

To conclude this part, we will apply the comparison theorem to prove the strong Feller property under Condition \ref{con1}.

\bgproposition\label{StrongFeller} Suppose that Condition~$\ref{con1}$ is satisfied. Then $(P_t)_{t\ge 0}$ has the strong Feller property.
\edproposition
\proof
Let $\{Y^{(1)}_t: t\ge 0\}$ and $\{Y^{(2)}_t: t\ge 0\}$ be two positive solutions to \eqref{sde1} with initial values $y_1>y_2\ge 0$ respectively, and let $Z_t=Y^{(1)}_t-Y^{(2)}_t$ for any $t>0$. It is clear that
$\{Z_t: t\ge 0\}$
takes value on $\R_+$ and solves the following equation
\beqlb\label{sdeZ}
Z_t\ar=\ar y_1-y_2 - \int_0^t bZ_{s}\,\d s + \sqrt{2c}\int_0^t\int_0^{Z_{s}} W_1(\d s,\d u) - \int_0^t [g(X_{s})-g(Y_{s})]\,\d s \cr
\ar\ar
+ \int_0^t\int_{0}^1\int_0^{Z_{s-}} z\,\tilde{M}_1(\d s,\d z,\d u) + \int_0^t\int_{1}^\infty\int_0^{Z_{s-}} z\,M_1(\d s,\d z,\d u),
\eeqlb
where $W_1(\d s,\d u):= W(\d s,Y^{(2)}_{s-}+\d u)$ is a Gaussian white noise on $(0,\infty)^2$ with intensity $\d s\d u$, and
$$M_1(\d s,\d z,\d u):= M(\d s,\d z,Y^{(2)}_{s-}+\d u)$$
is a Poisson random measure on $(0,\infty)^2$ with intensity $\d sm(\d z)\d u$.

We introduce a new process $(Z_t')_{t\ge0}$ by
\beqnn
Z'_t\ar=\ar y_1-y_2 - \int_0^t bZ'_{s}\,\d s + \sqrt{2c}\int_0^t\int_0^{Z'_{s}} W_1(\d s,\d u)\cr
\ar\ar
+ \int_0^t\int_{0}^1\int_0^{Z'_{s-}} z\,\tilde{M}_1(\d s,\d z,\d u) + \int_0^t\int_{1}^\infty\int_0^{Z'_{s-}} z\,M_1(\d s,\d z,\d u).
\eeqnn
Then  $(Z_t)_{t\ge0}$ is a CB-process with branching mechanism $\itPsi$ and its transition kernel $Q_t(x,\d y)$ satisfies \eqref{CB-semigr}.
It is clear that $(Z_t)_{t\ge0}$ can be controlled  by $(Z'_t)_{t\ge0}$; that is,
\[\mbb{P}(Z_t\le Z'_t)=1.\]

Note that $Z_t=0$ for $t\ge T$, where $T:=\inf\{t>0: Z_t=0\}$. Then, for any $f\in b\mcr{B}[0,\infty),$
\beqnn
|P_tf(y_1)-P_tf(y_2)|\ar\le\ar
\left|\mbb{E}\left[(f(Y^{(1)}_t)-f(Y^{(2)}_t))\I_{\{T>t\}}\right]\right|\\
\ar\le\ar 2\|f\|_\infty Q_t(y_1-y_2, (0, \infty)).
\eeqnn
It follows from \cite[Theorem 3.5]{Li11} that if Condition~$\ref{con1}$ is satisfied, then
\[Q_t(y_1-y_2, (0, \infty)) = 1-\e^{-(y_1-y_2)\bar{v}_t}\to 0 \quad\text{ as}\quad y_1-y_2\to 0,\]
where $\bar{v}_t<\infty$ for all $t>0.$ As a consequence, we know the function $x\mapsto P_tf(x)$ is continuous on $[0, \infty)$.
\qed

By directly applying Propositions $\ref{p1}$, $\ref{p2}$ and $\ref{StrongFeller}$, we can prove Propositions \ref{pc1}.

\subsection{Trajectory Feller property}
\medskip\noindent\textit{Proof of Proposition~$\ref{pc2}$.~}
As in the proof of Proposition \ref{StrongFeller}, let $(Y^{(1)}_t)_{t\ge 0}$ and $(Y^{(2)}_t)_{t\ge 0}$ be two positive solutions to \eqref{sde1} with initial values $y_1>y_2\ge 0$ respectively, and let $Z_t=Y^{(1)}_t-Y^{(2)}_t$ for any $t\ge0$.
Then, $(Z_t)_{t\ge 0}$ takes values on $\R_+$ and satisfies \eqref{sdeZ}.

Applying It\^o's formula in \eqref{sdeZ} and the facts that $x\e^{-x}\le 1-\e^{-x}$ for all $x\ge0$ and the function $x\mapsto g(x)$ is increasing, we get
\beqlb\label{Feller1}
1-\e^{-Z_t}\ar=\ar 1-\e^{-(y_1-y_2)}-\itPsi(1)\int^t_0 Z_{s}  e^{-Z_{s}} \d s \ccr
\ar\ar\, -\int^t_0 [g(Y^{(1)}_{s})-g(Y^{(2)}_{s})] Z_{s} \e^{-Z_{s}} \d s+M_t\ccr
\ar\le\ar
1-\e^{-(y_1-y_2)}+|\itPsi(1)|\int^t_0 (1-\e^{-Z_{s}}) \d s+M_t,
\eeqlb
where $\itPsi$ is defined in (\ref{bran-mech}), and
$$
M_t=\sqrt{2c}\int^t_0 \int_0^{Z_{s-}} \e^{-Z_{s-}} W_1(\d s, \d u)+\int^t_0\int^\infty_0\int^{Z_{s-}}_0 \e^{- Z_{s-}}[1-\e^{-z}] \tilde{M}_1(\d s, \d u, \d z)
$$
is a martingale satisfying
\[\mbb{E}[M^2_t]=\left[2c+\int_0^\infty (1-\e^{-z})^2 m(\d z)\right]\mbb{E}\left[ \int^t_0 Z_{s}\e^{-2 Z_s}\d s \right].\]

Taking  expectations  on both sides of (\ref{Feller1}), we obtain
\beqnn
\mbb{E}\left[1-\e^{-Z_t}\right]\le 1-\e^{-(y_1-y_2)}+|\itPsi(1)|\int^t_0 \mbb{E}[1-\e^{- Z_{s}}] \d s.
\eeqnn
Then, it follows from Gronwall's inequality that for any $t\ge0$,
\beqlb\label{Feller2}
\mbb{E}\left[1-\e^{-Z_t}\right]\le (1-\e^{-(y_1-y_2)})\e^{|\itPsi(1)|t}.
\eeqlb

On the other hand, taking supremes and expectations  on both sides of (\ref{Feller1}), we obtain
\beqlb\label{Feller5}
\mbb{E}\left[\sup_{0\le t\le T} \left(1-\e^{-Z_t}\right)\right]\ar\le\ar  1-\e^{-(y_1-y_2)}+\mbb{E}\left[\sup_{0\le t\le T}M_t\right]\ccr
\ar\ar +|\itPsi(1)|\int^T_0 \mbb{E}[1-\e^{-Z_{s}}] \d s.
\eeqlb
It follows from Jensen's inequality and  Burkholder-Davis-Gundy's inequality that
\beqlb\label{Feller4}
\left\{\mbb{E}\left[\sup_{0\le t\le T}M_t\right]\right\}^{2}
\ar\le\ar
\mbb{E}\left[\sup_{0\le t\le T}M^2_t\right]
\le
4\mbb{E}\left[M^2_T\right]\cr
\ar\le\ar
4\left[2c+\int_0^\infty (1-\e^{- z})^2 \mu(\d z)\right]\mbb{E}\left[ \int^T_0 Z_{s}\e^{-2 Z_s}\d s \right].
\eeqlb
Using \eqref{Feller2} and the fact that $x\e^{-2x}\le x\e^{-x}\le 1-\e^{-x}$ for all $x\ge 0$, we obtain
\beqlb\label{Feller3}
\mbb{E}\left[ \int^T_0 Z_{s}\e^{-2 Z_s}\d s\right]\le \mbb{E}\left[ \int^T_0 Z_{s}\e^{-Z_s}\d s \right]\le T\e^{|\itPsi(1)|T}(1-\e^{-(y_1-y_2)}).
\eeqlb

Substituting \eqref{Feller2}, \eqref{Feller4} and \eqref{Feller3} into \eqref{Feller5}, we know that for any $T>0$,
there exists a constant $C(T)>0$ such that
\[
\mbb{E}\left[\sup_{0\le t\le T} \left(1-\e^{-Z_t}\right)\right]\le C(T)(1-\e^{-(y_1-y_2)})^{1/2},
\]
which implies that in probability
\[\lim_{y_1-y_2\to0}\sup_{0\le t\le T}|Y^{(1)}_t-Y^{(2)}_t|=0.
\]
Noting that the Skorokhod distance is dominated by the supreme norm, we complete the proof.\qed

\subsection{Lyapunov condition}

\bglemma\label{pc3-r}
Suppose that there is an increasing and strictly positive function $g_0\in C^1[0,\infty)$ so that
\begin{itemize}
	\item[{\rm (i)}]  $$\int_{1}^\infty  \frac{1}{g_0(z)}\,\d z<\infty;$$
	\item [{\rm (ii)}] $$\limsup_{x\to \infty}\left[{ x\int_x^\infty \mu(\d  z})\int_x^\infty \frac{1}{g_0(z)}  \,\d z\right]<\infty$$ and
	$$\limsup_{x\to\infty}\frac{x\left(|b|+\int_ 1^x z\,\mu(\d z)\right)}{g_0(x) }<\infty;$$
	\item [{\rm (iii)}] $$ \lim_{x\to \infty} \frac{g(x)}{g_0(x)}=\infty.$$
\end{itemize}
Then {\bf (C3)} is satisfied.
\edlemma

\begin{proof} Define $f\in C^2_b[0, \infty)$ such that $f\ge1$,  $f'\ge 0$ and $f''\le 0$ on $\R_+$, and
	$$f(x)= 1+ \int_1^x \frac{1}{g_0(y)   }\,{\d y},\quad x\ge 1. $$
	It follows from (i) that $$1\le f(x) \le 1+ \int_1^\infty \frac{1}{g_0(y)  }\,{ \d y}=:C_0,\quad x\ge0.$$
	Then, applying the mean value theorem, we obtain that for any $x\ge 1$,
	\beqnn
	Lf(x) \ar=\ar  -[bx+g(x)]f^\prime(x)+cxf^{\prime\prime}(x)\cr
	\ar\ar
	+\,  x\int_0^\infty \big[f(x+z) - f(x) - zf^\prime(x)\I_{\{z\le 1\}}\big] \mu(\d z)\\
	\ar=\ar \left[ - bx-g(x)+x\int_1^x z\,\mu(\d z)\right] f^\prime(x)+cxf^{\prime\prime}(x)\cr
	\ar\ar
	+\,  x\int_0^x \big[f(x+z) - f(x) - zf^\prime(x) \big] \mu(\d z)\\
	\ar\ar +  x\int_x^\infty \big[f(x+z) - f(x)  \big] \mu(\d z)\\
	\ar\le\ar  \frac{- bx-g(x)+x\int_1^x z\,\mu(\d z)}{g_0(x)  }  +x\int_x^\infty  \mu(\d z)\int_x^\infty \frac{1}{g_0(y)  }\,\d y.
	\eeqnn According to the inequality above, (ii) and (iii),  one can take $l\ge 1$ large enough such that for all $x\ge l$,
	$$Lf(x) \le -\frac{g(x)}{2  g_0(x) }=:-\frac{1}{2}h(x)\le -\frac{1}{2C_0}h(x)f(x),$$ where
	in the last inequality we used the fact that $f(x)\le C_0$ for all $x\ge0$.
	
	Note again that, by (iii), $h(x)\to\infty$ as $x\to\infty$.
	For $n=1,2,\dots$, let \[r_n=\frac{1}{2C_0}h^*(l+n-1),\quad K_n=[0, l+n-1],\quad b_n=2r_n+\sup_{x\in K_n} |Lf(x)|,\] where $h^*(r)=\inf_{x\ge r}h(x).$ Then, we see that {\bf (C3)} is satisfied  with $W(x)=f(x)$.\qed
\end{proof}

{\noindent\textit{Proof of Proposition~$\ref{pc3}$.~} {Let $g_0\in C^1[0,\infty)$ be an increasing and strictly positive function such that for all $x\ge 1$, }
	$$g_0(x)=\begin{cases}
		\displaystyle{x\varphi(x)},&\quad \alpha\in (1,2),\\
		\displaystyle{x\varphi(x)\log { (1+x)}},&\quad \alpha=1,\\
		\displaystyle{x^{2-\alpha}},&\quad \alpha\in(0,1).\end{cases}$$ Then the desired assertion  is a consequence of Lemma \ref{pc3-r}.}

\subsection{Weak Feller property}

\noindent\textit{Proof of Proposition~$\ref{pc4}$.~}
Let $f$ be a continuous bounded function with support contained in $(0, \infty)$. Set
\beqnn
g(x)=\left\{\begin{array}{ll}
	\displaystyle f(x), &\quad x>0, \cr
	\displaystyle0, &\quad x=0.
\end{array}\right.
\eeqnn
It is clear that $g\in b\mcr{B}[0,\infty)$ and
\beqnn
\I_{\{t<\tau^-_0\}}f(Y_t)\ar=\ar \I_{\{t<\tau^-_0\}}g(Y_t)\\
\ar=\ar
\I_{\{t<\tau^-_0\}}g(Y_t)+\I_{\{t\ge\tau^-_0\}}g(Y_t)\\
\ar=\ar
g(Y_t),
\eeqnn
where we used the fact that $Y_t=0$ for all $t\ge \tau^-_0.$
From the definition of the killed semigroup $(P_t^0)_{t\ge 0}$, we see that
\beqnn
P_t^{0}f(x)- P_t^{0}f(y)
\ar=\ar
\mbb{E}_x[\I_{\{t<\tau^-_0\}}f(Y_t)]-\mbb{E}_y[\I_{\{t<\tau^-_0\}}f(Y_t)]\\
\ar=\ar
P_tg(x)-P_tg(y).
\eeqnn
By the strong Feller property of $(P_t)_{t\ge 0}$, we know the right hand side tends to zero as $x-y\to 0$, which gives the desired assertion.\qed

\subsection{Irreducibility}\label{section2.4}
Let $(X_t)_{t\ge 0}$ be a CB-process with branching mechanism $\itPsi$; that is, $(X_t)_{t\ge 0}$ is the strong solution to the SDE \eqref{sde1} with the function $g$ replaced by $0$. Then by
Lamperti \cite{Lam67b}  there is a spectrally negative L\'evy process $(N_t)_{t\ge0}$
with log-Laplace exponent $\itPsi$, i.e., \[\mbb{E}_x[\e^{-\lambda N_t}]=\e^{-\lambda x+\itPsi(\lambda)t}\] such that
\[N_{\eta^{-1}(t)\wedge S^{-}_0}=X_{t},\quad t\ge0,\] where
\[S^-_0:=\inf\{t>0: N_t\le 0\}, \quad\eta(t):=\int^{t}_0 \frac{1}{N_s} \d s\quad\mbox{and}\quad \eta^{-1}(t):=\inf\{s>0: \eta(s)>t\}.\]
For any $z> \varepsilon\ge0$, let $N^\varepsilon_t:=N_{t\wedge S^-_\varepsilon}$ and $X^\varepsilon_t:=X_{t\wedge T^-_\varepsilon},$
and define stopping times
\[S^-_\varepsilon:=\inf\{t\ge0: N_t\le \varepsilon\}=\inf\{t\ge0: N^\varepsilon_t\le \varepsilon\},\]
\[T^-_\varepsilon:=\inf\{t\ge0: X_t\le \varepsilon\}=\inf\{t\ge0: X^\varepsilon_t\le \varepsilon\},\]
\[S^{\varepsilon, +}_z:=\inf\{t\ge0: N^\varepsilon_t\ge z\}\quad\mbox{and}\quad T^{\varepsilon, +}_z:=\inf\{t\ge0: X^\varepsilon_t\ge z\}.\]
Then, from the random time change, we have for all $0\le \varepsilon<z$,
\beqlb\label{time change}
T^{-}_\varepsilon=\int_0^{S^{-}_\varepsilon}\frac{1}{N^\varepsilon_t}\d t
\quad\mbox{and}\quad
T^{\varepsilon, +}_z=\int_0^{S^{\varepsilon,+}_z} \frac{1}{N^\varepsilon_t}\d t.
\eeqlb

\bglemma\label{l7}
Suppose that Condition $\ref{nontri}$ is satisfied. Then for all $t>0$ and $z,x>\varepsilon\ge0$, we have
\[
\mbb{P}_x(S^{-}_\varepsilon<t)>0\quad\mbox{and}\quad \mbb{P}_x(S^{\varepsilon, +}_z<t)>0.
\]
\edlemma
\proof It follows from Sato \cite[Theorem 24.10 (i)]{Sato99} that  under Condition \ref{nontri} for any $t>0$, the support of $N_t$ is $(-\infty, +\infty)$, and hence the first assertion holds.

We next turn to the proof of the second assertion. Let $S^+_z:=\inf\{t\ge0: N_t\ge z\}.$ It is sufficient to prove that for any $t$ small enough and any $0\le\varepsilon<x<z$, $\mbb{P}_x(S^{+}_z<t<S^{-}_\varepsilon)>0.$

If $(N_t)_{t\ge0}$ has a Brownian motion component $(B_t)_{t\ge0}$, then $(B_t)_{t\ge0}$ is independent of $(N_t-B_t)_{t\ge0}$. Without loss of generality we set $B_0=0$. Since $(N_t-B_t)_{t\ge 0}$ is a L\'evy process and has the same initial value as $(N_t)_{t\ge0}$, applying the right-continuity of $(N_t-B_t)_{t\ge 0}$, we know $\mbb{P}_x(\lim_{t\to 0}(N_t-B_t)=x)=1$. Then, for small enough $t>0$,
\beqnn
\begin{split}
	\mbb{P}_x(S^{+}_z<t<S^{-}_\varepsilon)\ge &\mbb{P}_x\left(\inf_{0\le s\le t} (N_s-B_s)>(x+\varepsilon)/2\right)\\
	&\times\mbb{P}\left(\sup_{0\le s\le t} B_s>z,\inf_{0\le s\le t} B_s>-(x-\varepsilon)/2\right)>0.\end{split}
\eeqnn

We next prove the case that $(N_t)_{t\ge0}$ is a pure jump process that can be expressed as
\[N_t=x-bt+\int^t_0\int^1_0 u \tilde{N}(\d s, \d u)+\int^t_0\int^\infty_1 u N(\d s, \d u),\]
where $N(\d s, \d u)$ denotes a Poisson random measure with intensity $\d s \mu(\d u)$.
We introduce a new process
\[N'_t=x-\left[b+\int^1_\kappa u \mu(\d u)\right]t+\int^\infty_0\int^\kappa_0 u \tilde{N}(\d s, \d u),\]
with $\kappa\in (0,1]$ such that $\mu[\kappa, \infty)>0$. It is clear that
\[N_t=N'_t+\int^\infty_0\int^\infty_\kappa u N(\d s, \d u).\]
Let $S^{\prime -}_\varepsilon=\inf\{t\ge0: N'_t\le \varepsilon\}$. Then the stopping time $S^{\prime -}_\varepsilon$ is independent of the process $(\int^t_0\int^\infty_\kappa u N(\d s, \d u))_{t\ge0}$, and for small $t>0$, $\mbb{P}_x(S^{\prime -}_\varepsilon>t)>0$, because of the right continuity of $t\mapsto N'_t$. Since the number of jumps associated with the Poisson integral  $\int^t_0\int^\infty_\kappa u N(\d s, \d u)$
follows the Poisson distribution, we have \[\mbb{P}\left(\int^t_0\int^\infty_\kappa u N(\d s, \d u)>z\right)>0,\quad z,t>0.\]
Therefore, for $t>0$ small enough,
\beqnn\mbb{P}_x(S^{+}_z<t<S^{-}_\varepsilon)\ar\ge\ar \mbb{P}_x(S^{+}_z<t<S'^{-}_\varepsilon)\\
\ar\ge\ar \mbb{P}\left(\int^t_0\int^\infty_\kappa u N(\d s, \d u)> z+\varepsilon\right)\mbb{P}_x(S'^{-}_\varepsilon>t)\\
\ar>\ar0.\eeqnn  The proof is complete. \qed
\bglemma\label{l6}
Suppose that Condition $\ref{nontri}$ is satisfied. Then for all $t>0$ and $z, x>\varepsilon>0$, we have
\[
\mbb{P}_x(T^{-}_\varepsilon<t)>0\quad\mbox{and}\quad \mbb{P}_x(T^{\varepsilon, +}_z<t)>0.
\]
\edlemma
\proof
It follows from \eqref{time change} that
\[T^{-}_\varepsilon=\int_0^{S^{-}_\varepsilon}\frac{1}{N_t}\d t\le \varepsilon^{-1}S^{-}_\varepsilon\]
and
\[T^{\varepsilon, +}_z=\int_0^{S^{\varepsilon, +}_z}\frac{1}{N_t}\d t\le \varepsilon^{-1}S^{\varepsilon, +}_z.\]
Applying Lemma \ref{l7}, we know for all $t>0$,
\[
\mbb{P}_x(T^{-}_\varepsilon<t)\ge \mbb{P}_x(S^{-}_\varepsilon<\varepsilon t)>0\quad\mbox{and}\quad \mbb{P}_x(T^{\varepsilon, +}_z<t)\ge \mbb{P}_x(S^{\varepsilon, +}_z<\varepsilon t)>0.
\] The desired assertion follows.\qed

{Recall that $(Y_t)_{t\ge0}$ is the CBC-process and \[\tau^-_z:=\inf\{t>0: Y_t\le z\},\qquad z\ge 0.\] Now, we let
\[
\sigma_z:=\inf\{t\ge0: Y_t=z\}, \qquad z\ge 0.
	\]}
We have the following result.
\bglemma\label{l4}
Suppose that Condition $\ref{nontri}$ is satisfied. Then for any $t,x,z>0$, we have
\[\mbb{P}_x(\sigma_z<t)>0.\]
\edlemma
\proof
Choose a constant $\varepsilon\in (0, x\wedge z)$. We first consider the case that $x>z>0$.
Since the CBC-process $(Y_t)_{t\ge0}$ has no negative jump, by a comparison argument and Lemma \ref{l6}, we know that $$\mbb{P}_x(\sigma_z< t)\ge \mbb{P}_x(\tau^-_\varepsilon< t)\ge \mbb{P}_x(T^-_\varepsilon< t)>0.$$

We next consider the situation that $0<x<z $. Since the function $g$ is locally bounded, there is a constant $h:=h(\varepsilon,z)>0$ such that $g(y)<hy$ for all $y\in [\varepsilon, z].$ Let $(X^h_t)_{t\ge 0}$ be the strong solution to \eqref{sde1} with $g(x)$ replaced by $hx$. Then $(X^h_t)_{t\ge 0}$ is a CB-process with the branching mechanism
\[\itPsi^h(\lambda)= (b+h)\lambda + c\lambda^2 + \int_0^\infty \big(\e^{-\lambda z}-1+\lambda z\I_{\{z\le 1\}}\big)\mu(\d z), \quad \lambda\ge 0.\]

Let $(Y^\varepsilon_t)_{t\ge0}:=(Y_{t\wedge \tau^-_\varepsilon})_{t\ge0}$ be the process $(Y_t)_{t\ge0}$ stopped at $\varepsilon$.
Define stopping times
\[\sigma^\varepsilon_z:=\inf\{t\ge0: Y^\varepsilon_t=z\},\quad \tau^{\varepsilon, +}_z:=\inf\{t>0: Y^\varepsilon_t\ge z\}\]
and
\[ T^{h,\varepsilon,+}_z:=\inf\{t>0: X^{h, \varepsilon}_t\ge z\}.\]
Then by a comparison result and Lemma \ref{l6} we have for any $t>0$ and $y>\varepsilon$,
\beqlb\label{ineq20}\mbb{P}_x(\tau^{\varepsilon, +}_z< t)\ge\mbb{P}_x(T^{h,\varepsilon, +}_z< t)>0
\eeqlb
and 
\beqlb\label{ineq21}\mbb{P}_{y}(\tau^{\varepsilon,-}_\varepsilon\le t/2)= \mbb{P}_{y}(\tau^-_\varepsilon\le t/2)\ge\mbb{P}_{y}(T^-_\varepsilon\le t/2)>0.\eeqlb
Here for any $y>y_0\ge \varepsilon$,
$$\tau^{\varepsilon, -}_{y_0}=\inf\{t>0: Y^\varepsilon_t\le y_0\}.$$

On the other hand,
it is clear that
\beqlb\label{ineq22}
\mbb{P}_x(\sigma_z<t)\ge\mbb{P}_x(\sigma_z<t\wedge \tau^-_\varepsilon)=\mbb{P}_x(\sigma^\varepsilon_z<t).
\eeqlb
Suppose that the process $(Y^\varepsilon_t)_{t\ge0}$  first upcrosses level $z$ during the  time interval $(0, t/2)$, and then hits $\varepsilon$ from above during $[\tau^{\varepsilon, +}_z, \tau^{\varepsilon, +}_z+t/2)$. Then it  has to reach $z$ before time $t$,  since the process $(Y^\varepsilon_t)_{t\ge0}$ has no negative jump. In other words, we have
\beqlb\label{ineq23}
\mbb{P}_x(\sigma^\varepsilon_z<t) \ge \mbb{P}_x(\tau^{\varepsilon, +}_z<t/2, \tau^{\varepsilon,-}_\varepsilon\circ \theta_{\tau^{\varepsilon, +}_z}<t/2)
\eeqlb
where $\theta$ denotes the shift operator. By the strong Markov property, we know the right-hand side can be written by
\beqlb\label{ineq24}
\mbb{P}_x(\tau^{\varepsilon, +}_z<t/2,  \tau^{\varepsilon,-}_\varepsilon\circ\theta_{\tau^{\varepsilon, +}_z}<t/2)
\ar=\ar
\mbb{E}_x\left[\I_{\{\tau^{\varepsilon, +}_z<t/2\}}\mbb{P}_{Y^\varepsilon_{\tau^{\varepsilon, +}_z}}(\tau^{\varepsilon,-}_\varepsilon< t/2) \right]
.\eeqlb

Furthermore, by \eqref{ineq21},
$$\mbb{P}_{Y^\varepsilon_{\tau^{\varepsilon, +}_z}}(\tau^{\varepsilon,-}_\varepsilon\le t/2)
>0$$ almost surely, which together with \eqref{ineq20} gives
\beqlb\label{ineq25}\mbb{E}_x\left[\I_{\{\tau^{\varepsilon, +}_z<t/2\}}\mbb{P}_{Y^\varepsilon_{\tau^{\varepsilon, +}_z}}(\tau^{\varepsilon,-}_\varepsilon< t/2) \right]
>0.\eeqlb
Finally, putting \eqref{ineq22}, \eqref{ineq23}, \eqref{ineq24} and \eqref{ineq25} together yields the desired result.\qed

\medskip\noindent\textit{Proof of Proposition~$\ref{pc5}$.~}
Fix $z\in O$ and $t>0$. By Lemma \ref{Lemma-uniform}, we can find a small enough  $s\in(0, t)$ such that $P_u(z, O)>0$ for all $u\le s$. It follows from Lemma \ref{l4} that for any $y>0$, $\mbb{P}_y(\sigma_z<s)>0$. Applying the strong Markov property of $(Y_t)_{t\ge0}$ at time $\sigma_z$, we know that for any $y>0$,
$$P_{s}(y, O)\ge\mbb{E}_y\left[\I_{\{\sigma_z<s\}}P_{s-\sigma_z}(z, O)\right]>0.$$
On the other hand, there exist a positive integer $N$ and a real number $u\in (0, (t-s)\wedge s)$ such that $t-s=Nu$. Then
$$
P_{t-s}(x, (0, \infty))=\int_{(0, \infty)^N} P_{u}(x, \d y_1)P_{u}(y_1, \d y_2)\cdots P_{u}(y_{N-1}, \d y_N)>0.
$$
Finally,
$$P_t^0(x, O)=P_{t}(x, O)=\int_{(0, \infty)} P_{s}(y, O) P_{t-s}(x, \d y)>0.$$

We next turn to prove the second assertion. Suppose that Condition \ref{con1} holds. Then by a comparison argument and \cite[Theorem 1]{Grey}  we know that for all $x, t>0$, $$ \mbb{P}_x(\tau^-_0< t)\ge \mbb{P}_x(T^-_0< t)>0.$$
Suppose that Condition \ref{con2} holds.  Fix $x, t>0$, according to \eqref{neweq1}, there is a $z\in (0, x)$ such that
	\[\mbb{P}_z (\tau^-_0<t/2)>0.\]
	Then by the strong Markov property and Lemma \ref{l4}, we obtain
	\[\mbb{P}_x(\tau^-_0<t)\ge \mbb{P}_x (\sigma_z<t/2)\mbb{P}_z (\tau^-_0<t/2)>0.\]
The proof is completed.
\qed

\noindent\textbf{Acknowledgements.~} This research is supported by the National Key R{\&}D Program of China (Nos.~2022YFA1000033), the National Natural Science Foundation of China (Nos. ~11831014, 12071076, 12225104 and 12271029), the Education and Research Support Program for Fujian Provincial Agencies and Natural Science ans Engineering Research Council of Canada (RGPIN-2021-04100).

\end{document}